\UseRawInputEncoding
\documentclass[a4paper]{article}
\usepackage{amssymb}
\usepackage{amsmath}
\usepackage{amsthm}
\usepackage{color}
\usepackage{geometry}
\usepackage{graphicx}
\usepackage{epsf}
\usepackage{subfigure}
\usepackage{multirow}
\newtheorem{thm}{Theorem}[section]
\newtheorem{lam}[thm]{Lemma}
\newtheorem{Def}[thm]{Definition}
\newtheorem{Rem}[thm]{Remark}
\makeatletter
\newcommand{\rmnum}[1]{\romannumeral #1}
\makeatother

\title{Existence and convergence of solutions for nonlinear elliptic systems on graphs}

\begin{document}

\date{2021.11.20}

\author{Jinyan Xu$^{1}$\ \ \ Liang Zhao$^{1}$\\
\it\small  $^{1}$  School of Mathematical Sciences, Key Laboratory of Mathematics and Complex Systems of MOE, \\
\it\small Beijing Normal University, Beijing 100875, China\\
}

\maketitle

\begin{abstract}
We consider a kind of nonlinear systems on a locally finite graphs $G=(V,E)$. We prove
via the mountain pass theorem that this kind of systems has a nontrivial ground state
solution which depends on the parameter $\lambda$ with some suitable assumptions on the
potentials. Moreover, we pay attention to the concentration behavior of these
solutions and prove that, as $\lambda \to \infty$, these solutions converge to a ground
state solution of a corresponding Dirichlet problem. Finally, we also provide some
numerical experiments to illustrate our results.

\textbf{Keywords:} nonlinear elliptic system, locally finite graph, ground state solution
\end{abstract}

\section{Introduction}

Recently, the analytic aspect of different partial differential equations on graphs has attracted much attention. For example, a variety of fundamental problems of heat equations on graphs, such as the heat kernel \cite{Horn}, the existence, uniqueness \cite{Huang, Lin1} and blow-up behavior \cite{Lin2} of solutions are considered by different authors.
There is also research about the Fokker-Planck and Schr\"odinger equations on graphs \cite{Chow1, Chow2}, which are related to the discrete optimal
transport theory. Grigoryan, Lin and Yang \cite{Gri1,Gri2,Gri3} studied several
nonlinear elliptic equations on graphs and they pointed out that the required Sobolev
spaces on a finite graph is pre-compact, which makes it possible to use the variational method to obtain the existence of solutions. The second author cooperates with others \cite{ZhangZhao, Han} by using the Nehari manifold to prove that on a locally finite graph, the nonlinear Schr\"odinger equation has a nontrivial ground state solution under suitable conditions, and the limit of the solution is limited to a potential well.

A single nonlinear partial differential equation defined on Euclidean space as
$$-\Delta u+b\left (x\right )u=f\left (x,u\right ),$$
it has been extensively studied. For example, Rabinowitz\cite{Rab} proposed and
studied the standing wave solution of the Schr\"odinger equation. Li\cite{Li} and Bartsch and Willem\cite{BarWil} also studied the existence of solutions under different assumptions. In particular, Bartsch and Wang \cite{BarWang1,BarWang2} considered the potential $b\left (x\right )=\lambda a\left (x\right )+1$ where $a\left(x\right)$ satisfies certain assumptions and they proved that the equation has a ground state solution that depends on
the parameter $\lambda$, which will converge to a ground state solution of a Dirichlet
problem when $\lambda \to \infty$. For different assumptions on the potential
$a\left(x\right)$, we refer to \cite{FurMaiSil,BarTang} and the references therein.

On the other hand, systems of nonlinear partial differential equations, which are obviously generalizations of a single equation, is used to describe many phenomena in nature, such as nerve conduction and electromagnetic fields et al. Kinds of nonlinear systems defined on Euclidean space have been extensively studied. For the general gradient systems
\begin{displaymath}
\left\{
\begin{array}{ll}
-\Delta u+a\left (x\right )=F_{u}\left (x,u,v\right )
\\
-\Delta u+b\left (x\right )=F_{v}\left (x,u,v\right )
\\
u,v~~\in ~~W^{1,2}\left (\mathbb{R}^{N} \right )
\end{array}
\right.
\end{displaymath}
when $a\left (x\right )=b\left (x\right )\equiv 1$, the existence of weak solutions is proved in \cite{Kri}. In \cite{Cos}, Costa required the potentials to be continuous and coercive to overcome the lack of compactness and proved the existence of weak solutions. For the nonlinear gradient systems, we can also refer to \cite{Alv1,Fig, LiuLiu, Ou, ZhangZhang}.

In \cite{Alv2}, Alves considered the following system with a parameter $\epsilon$
\begin{displaymath}
 \left\{
\begin{array}{ll}
-\epsilon ^{2}\Delta u+a\left (x\right )=F_{u}\left (x,u,v\right )
& in ~~ \mathbb{R}^{N}
\\
-\epsilon ^{2}\Delta u+b\left (x\right )=F_{v}\left (x,u,v\right )
& in ~~ \mathbb{R}^{N}
\\
u\left (x\right ),v\left (x\right ) \to 0 ~~as~~\left | x \right | \to \infty
\\
u,v>0
& in ~~ \mathbb{R}^{N}
\end{array}
\right.
\end{displaymath}
Using the mountain pass theorem, the author showed the existence results and the concentration behavior of these solution as $\epsilon \to 0$. There has been a great deal of interest in studying the convergence of solutions for nonlinear systems with different assumptions, we refer \cite{FigFur,Lv,Ou} and the references therein.

We are concerned in this article with the following nonlinear system on a locally finite and connected graph $G=\left(V,E\right)$.
\begin{equation}\label{system}
\left\{\begin{array}{ll}
- \Delta u+\left ( \lambda a\left ( x \right ) +1 \right ) u=\frac{\alpha }{\alpha +\beta}
\left | u \right | ^{\alpha  -2}u\left | v \right |  ^{\beta }
\qquad
&\textrm{in}  ~V
 \\
-\Delta v+\left ( \lambda b\left ( x \right ) +1 \right ) v=\frac{\beta}{\alpha +\beta}
\left | u \right | ^{\alpha } \left | v \right |  ^{\beta -2}v
\qquad
&\textrm{in}  ~V
\end{array}\right.
\end{equation}
where $\alpha $ and
$\beta$ are positive constants such that $\alpha ,\beta >1$, $\lambda >0$ is a
parameter, and $a(x),b(x):V \to \mathbb{R}$ are given functions satisfying the following two conditions:
\begin{description}
\item[$\left(A_{1}\right)$]   $a\left(x\right)\geqslant 0, b\left(x\right)\geqslant 0$
for all $x\in V$, $\Omega _{a}:=\left \{ x\in V:a\left ( x \right ) =0 \right \},
\Omega _{b}:=\left \{ x\in V:b\left ( x \right ) =0 \right \} $ and
$\Omega _{a}\cap \Omega _{b}$ are all non-empty bounded domains in $V$.
\item[$(A_{2})$]  There exists a vertex $x_{0}\in V$ such that $a\left(x\right)\to
+\infty $ and $b\left(x\right)\to +\infty $ as $d\left(x,x_{0}\right)\to + \infty$.
\end{description}

Motivated by \cite{Han,ZhangZhao}, we aim to study the existence of a nontrivial ground state solution of (\ref{system}) with a fixed positive parameter $\lambda$. Here a ground state solution means it has the least energy among all nontrivial solutions. Moreover, we also study the convergence behavior of the solutions as $\lambda \to \infty$.

To describe this problem in details, we first introduce some concepts and assumptions. Let $G=\left(V,E\right)$ be a graph, where $V$ denotes the set of vertices and $E$ denotes the set of edges, and we write $x\sim y$ if $x$ is connected to $y$, i.e.  $xy \in E$. A graph $G$ is called \emph{locally finite} if each vertex has a finite number of edges. A graph is called \emph{connected} if any two vertices $x$ and $y$ can be connected via finite edges. The \emph{graph distance} $d\left ( x,y \right ) $ between any two distinct vertices $x,y$ is the minimal number of edges which connect these two vertices. We use $w_{xy}>0$ to denote the weight of an edge $xy \in E$ and we call it a \emph{symmetric weight} on $G$, if $w_{xy}=w_{yx}$ for any $\left ( x,y \right ) \in E$. The \emph{measure} $\mu :V\to \mathbb{R}^{+} $ on the graph is a finite positive function on $G$. We call it a uniformly positive measure if there exists a constants $\mu _{\min}>0$ such that $\mu \left ( x \right ) \geqslant \mu _{\min}$ for all $x\in V$.
If the distance $d\left ( x,y \right ) $ is uniformly bounded from above for any $x,y \in \Omega$, we call $\Omega$ a bounded domain in $V$. The boundary of $\Omega$ is defined as
\begin{displaymath}
\partial \Omega :=\left \{ y \notin \Omega :\exists x\in \Omega ~~such ~that~ xy \in E \right \}
\end{displaymath}
and the interior of $\Omega$ is denoted by $\Omega^{0}$. Moreover, we denote $\overline{\Omega}=\Omega \cup \partial \Omega$.

For any function $u:V\to \mathbb{R}$, the $\mu$-Laplacian of $u$ at $x$ is defined by:
\begin{displaymath}
\Delta u( x ) :=\frac{1}{\mu ( x ) } \sum_{y\sim x}^{} w_{xy}( u(y)-u(x)  ).
\end{displaymath}

The gradient form of two functions $u$ and $v$ on the graph is defined by:
\begin{displaymath}
\Gamma (u,v)(x):=\frac{1}{2\mu(x)}\sum_{y\sim x}^{}  w_{xy}(u(y)-u(x))(v(y)-v(x)).
\end{displaymath}

In particular, we use $\Gamma(u)$ to denote $\Gamma(u,u)$ and the length of the gradient for $u$ at $x$ is
\begin{displaymath}
\left | \nabla u \right | (x):=\sqrt{\Gamma(u)(x)}=\left(\frac{1}{2\mu(x)}\sum_{y\sim x}^{}w_{xy} \big (u(y)-u(x)\big )^{2} \right )^{1/2}.
\end{displaymath}

The integral of $u$ over $V$ is defined by:
\begin{displaymath}
\int_{V}^{}ud\mu=\sum_{x\in V}^{}  \mu(x)u(x).
\end{displaymath}

The set of functions with compact support is $C_{c}\left(V\right):=\left\{u:V\to \mathbb{R} ~|~\left \{ x\in V:u\left ( x \right )\ne0 \right \} \emph{\emph{is of finite cardinality}} \right\}$. In \cite{Han}, they have proved that $W^{1,2}\left(V\right)$ under the norm
\begin{displaymath}
\left \| u \right \| _{W^{1,2}\left ( V \right ) }=\left ( \int_{V}^{} \left ( \left | \nabla u \right |^{2}+u^{2}  \right ) d\mu  \right )^{1/2}
\end{displaymath}
is the completion of $C_{c}\left(V\right)$. Let us consider the space $H=W^{1,2}\left ( V \right ) \times W^{1,2}\left ( V \right ) =\left \{ \left ( u,v \right ) ~|~ u,v\in W^{1,2}\left ( V \right )   \right \} $.
It is easy to prove that $H$ is a Hilbert space with the inner product
\begin{displaymath}
\left \langle \left ( u,v \right ) ,\left ( \xi ,\eta  \right )  \right \rangle _{H}=\int_{V}^{} \left ( \Gamma \left ( u,\xi \right ) +u\xi +\Gamma \left ( v,\eta \right ) +v\eta  \right )d\mu
\end{displaymath}
for any $\left ( u,v \right ) ,\left ( \xi,\eta \right ) \in H$, and it is the completion of $C_{c}\left(V\right) \times C_{c}\left(V\right)$ under the norm
\begin{displaymath}
\left \| \left ( u,v \right )  \right \| _{H}^{2} =\int_{V}^{} \left ( \left | \nabla u \right | ^{2}+\left | \nabla v \right | ^{2}+u^{2}+v^{2} \right )d\mu
\end{displaymath}

Next, in order to apply the variational setting, for any fixed $\lambda>0$, we introduce the following subspace $H_{\lambda }$ of $H$:
\begin{displaymath}
H_{\lambda}=\left \{ \left ( u,v \right ) \in H ~|~ \int_{V}^{} \left (\lambda a\left ( x \right ) u^{2}+\lambda b\left ( x \right ) v^{2} \right)d\mu <\infty  \right \} .
\end{displaymath}
where $a\left ( x \right ) $ and $b\left ( x \right ) $ satisfy $\left( A_{1} \right)$ and $\left( A_{2} \right)$. Obviously, $H_{\lambda}$ is also a Hilbert space with its inner product and norm given by
\begin{displaymath}
\left \langle \left ( u,v \right ) ,\left ( \xi ,\eta  \right )  \right \rangle _{H_{\lambda}}=\int_{V}^{} \left ( \Gamma \left ( u,\xi \right ) + \left ( \lambda a\left ( x \right ) +1 \right ) u\xi +\Gamma \left ( v,\eta \right ) +\left ( \lambda b\left ( x \right ) +1 \right ) v\eta  \right )d\mu
\end{displaymath}
\begin{displaymath}
\left \| \left ( u,v \right )  \right \| _{H_{\lambda}}^{2} =\int_{V}^{} \left ( \left | \nabla u \right | ^{2}+\left | \nabla v \right | ^{2}+\left ( \lambda a\left ( x \right ) +1 \right )  u^{2}+\left ( \lambda b\left ( x \right ) +1 \right )  v^{2} \right )d\mu
\end{displaymath}
for any $\left ( u,v \right ) ,\left ( \xi ,\eta  \right ) \in H_{\lambda }$.

The energy functional $J_{\lambda}:H_{\lambda}\to \mathbb{R}$ associated with the system (\ref{system}) is defined by:
\begin{align*}
J_{\lambda}\left ( u,v \right ) :=&\frac{1}{2} \int_{V}^{}  \left ( \left | \nabla u \right | ^{2}+\left | \nabla v \right | ^{2}+\left ( \lambda a\left ( x \right ) +1 \right ) u^{2}+\left ( \lambda b\left ( x \right ) +1 \right ) v^{2} \right )d\mu
\\
& - \frac{1}{\alpha +\beta}\int_{V}^{}\left | u \right |^{\alpha } \left | v \right |^{\beta } d\mu
\\
=&\frac{1}{2} \left \| \left ( u,v \right )  \right \| _{H_{\lambda }}^{2}-\frac{1}{\alpha +\beta}\int_{V}^{}\left | u \right |^{\alpha } \left | v \right |^{\beta } d\mu
\end{align*}

In view of the hypotheses $\left( A_{1} \right)$ and $\left( A_{2} \right)$, the functional $J_{\lambda}$ is well defined and of class $C^{1}$. In order to obtain the critical points of $J_{\lambda}$ which are weak solutions of (\ref{system}), we will verify the Palais-Smale condition and apply the Mountain Pass Theorem to guarantee critical points of the functional. Our main results can be formulated as follows.

\begin{thm}\label{ssolution}
Let $G=\left(V,E\right)$ be a locally finite and connected graph with symmetric weight and uniformly positive measure. Assume $a\left ( x \right ) ,b\left ( x \right ) $ are functions satisfying $\left ( A_{1} \right ) $ and $\left ( A_{2} \right )$. Then for any positive constant $\lambda > 0$, there exists a ground state solution $\left(u_{\lambda},v_{\lambda}\right)$ of the system (1).
\end{thm}

To study the behavior of $u_{\lambda}$ as $\lambda \to \infty$, we introduce the Dirichlet problem:
\begin{equation}\label{dirichlet}
\left \{
\begin{array}{cl}
-\Delta u + u = \frac{\alpha }{\alpha +\beta } \left | u \right | ^{\alpha -2}u\left | v \right | ^{\beta }
\qquad
&\textrm{in} ~\Omega _{a}
\\
-\Delta v + v = \frac{\beta }{\alpha +\beta } \left | u \right | ^{\alpha}u\left | v \right | ^{\beta -2}v
\qquad
&\textrm{in} ~\Omega _{b}
\\
u = 0
&\textrm{on} ~ \partial  \Omega _{a}
\\
v = 0
&\textrm{on} ~ \partial  \Omega _{b}
\end{array}
\right .
\end{equation}
It is suitable to study (\ref{dirichlet}) in the space $H_{\Omega} :=W_{0}^{1,2}\left ( \Omega_{a}  \right )\times W_{0}^{1,2}\left ( \Omega _{b} \right )$ where $W_{0}^{1,2}\left ( \Omega  \right )$ is the completion of $C_{c}\left ( \Omega  \right )$ under the norm
\begin{displaymath}
\left \|  u   \right \| _{W_{0}^{1,2}\left ( \Omega  \right )} =\left (\int_{\Omega \cup \partial \Omega }^{}  \left | \nabla u \right | ^{2}+\int _{\Omega }u^{2} d\mu\right )^{\frac{1}{2}  }
\end{displaymath}
where $C_{c}\left ( \Omega  \right )$ denotes the set of all functions $u:\Omega \to \mathbb{R}$ satisfying $supp ~ u \subset \Omega$ and $u = 0$ on $\partial \Omega$. The space $H_{\Omega}$ endowed with the inner product
\begin{displaymath}
\left \langle \left ( u,v \right ) ,\left ( \xi ,\eta  \right )  \right \rangle _{H_{\Omega}}=\int_{\overline{\Omega}_{a} \cup \overline{\Omega}_{b} }^{} \left ( \Gamma \left ( u,\xi \right ) +\Gamma \left ( v,\eta \right ) \right )d\mu +\int _{\Omega_{a}\cup\Omega_{b} } \left (u\xi +v\eta  \right )d\mu
\end{displaymath}
is a Hilbert space.The functional related to (\ref{dirichlet}) is
\begin{displaymath}
J_{\Omega}\left (u,v\right )=\frac{1}{2} \int _{\overline{\Omega}_{a} \cup \overline{\Omega}_{b}}\left (\left | \nabla u \right |^{2}+\left | \nabla v \right |^{2} \right )d\mu
+\frac{1}{2} \int _{\Omega_{a}\cup\Omega_{b} }\left (u^{2}+v^{2}\right )d\mu
-\frac{1}{\alpha +\beta }\int _{\Omega_{a}\cup\Omega_{b} }\left | u \right |  ^{\alpha }\left | v \right | ^{\beta }d \mu
\end{displaymath}

Similar to Theorem \ref{ssolution}, the system (\ref{dirichlet}) also has a ground state solution. Precisely, we have

\begin{thm}\label{dsolution}
Let $G = \left( V, E \right)$ to be a locally finite and connected group with symmetric weight and uniform positive measure. Suppose that $\Omega_{a}$, $\Omega_{b}$ and $\Omega _{a}\cap \Omega _{b}$ are non-empty and bounded domains in $V$. Then the system (\ref{dirichlet}) has a ground state solution $\left( u_\Omega,v_\Omega \right) \in H_{\Omega}$.
\end{thm}

Finally, as in the scale case \cite{Ou}, the next result show that the semilinear elliptic system (\ref{dirichlet}) may be seen as a limit problem for (\ref{system}) when $\lambda$ goes to infinity, where $\Omega_{a}$ and $\Omega_{b}$ are the potential wells of the system (\ref{system}). More precisely, we have the following theorem.

\begin{thm}\label{convergence}
Under the same assumptions as in Theorem \ref{ssolution} and 1.2, we have that , for any sequence $\lambda_{n} \to \infty$, up to a subsequence, the corresponding ground state solutions $\left( u_{\lambda_{n}},v_{\lambda_{n}}\right)$ of (\ref{system}) converge in $H$ to a ground state solution of (\ref{dirichlet}).
\end{thm}

The paper is organized in the following way. In section 2 we present some basic properties and known results which will be used throughout our work. We prove Theorem \ref{ssolution} and \ref{dsolution} in section 3 and section 4 is devoted to the proof of Theorem \ref{convergence}. In the final section, we give a numerical experiment on a finite graph with $22$ vertices to illustrate our conclusions.

\section{Preliminaries}

\subsection{Weak solution}

We say that $\left( u,v \right)$ is a \emph{solution} of the system (\ref{system}) if the two equations hold for all $x\in V$. To define the weak solution, we need formulas of integration by parts on graphs, which are also fundamental when we use methods from calculus of variations. The proofs of the next two lemmas can be found in \cite{Han} and we omit them here.

\begin{lam}
Suppose that $u \in W^{1,2}\left(V\right)$. Then for any $\xi \in C_{c}\left(V\right)$ we have
\begin{displaymath}
\int_{V}\nabla u\nabla \xi d\mu =\int_{V} \Gamma \left (u,\xi \right )d\mu =-\int_{V} \left ( \Delta u \right )\xi d\mu
\end{displaymath}
\end{lam}
\begin{lam}
Suppose that $u \in W_{0}^{1,2}\left( \Omega \right)$. Then for any $\xi \in C_{c}\left(\Omega \right)$ we have
\begin{displaymath}
\int_{\Omega \cup \partial \Omega }\nabla u\nabla \xi d\mu =\int_{\Omega \cup \partial \Omega} \Gamma \left (u,\xi \right )d\mu =-\int_{\Omega \cup \partial \Omega} \left ( \Delta u \right )\xi d\mu
\end{displaymath}
\end{lam}

Now we can define the weak solution of the system (\ref{system}).

\begin{Def}
Suppose $\left( u,v \right) \in H_{\lambda}$. If for any $\left(\xi ,\eta \right) \in H_{\lambda}$, there holds
\begin{align*}
&\int_{V}\left(\nabla u\nabla \xi +\nabla v\nabla \eta +\left (\lambda a\left (x\right )+1\right )u\xi +\left (\lambda b\left (x\right )+1\right )v\eta \right)d\mu
\\
=& \int _{V}\left(\frac{\alpha }{\alpha +\beta }\left |u\right |^{\alpha -2}u\left |v\right |^{\beta }\xi +\frac{\beta }{\alpha +\beta }\left |u\right |^{\alpha }\left |v\right |^{\beta -2} v\eta \right)d\mu,
\end{align*}
then $\left(u,v\right)$ is called a weak solution of (1).
\end{Def}

Similarly, the weak solution of the system (2) is defined as
\begin{Def}
Suppose $\left( u,v \right) \in H_{\Omega}$. If for any $\left(\xi ,\eta \right) \in H_{\Omega}$, there holds
\begin{displaymath}
\int_{\overline{\Omega}_{a} \cup \overline{\Omega}_{b}}\left(\nabla u\nabla \xi +\nabla v\nabla \eta \right)d\mu
+\int_{\Omega_{a}\cup\Omega_{b} }\left(u\xi +v\eta \right)d\mu
= \int _{\Omega_{a}\cup\Omega_{b}}\left(\frac{\alpha }{\alpha +\beta }\left |u\right |^{\alpha -2}u\left |v\right |^{\beta }\xi +\frac{\beta }{\alpha +\beta }\left |u\right |^{\alpha }\left |v\right |^{\beta -2} v\eta \right)d\mu
\end{displaymath}
then $\left( u,v \right) $ is called a weak solution of (2).
\end{Def}

\subsection{Sobolev embedding}
Next, we present some results about the compactness of the function spaces $H_{\lambda}$ and $H_{0}$.
\begin{lam}\label{embedding1}
Assume that $\alpha,\beta>1, \lambda>0$ and $a\left(x\right),b\left(x\right)$ satisfies $\left(A_{1}\right)$ and $\left(A_{2}\right)$. Then $H_{\lambda}$ is embedded continuously into $H$ and $L^{q}\left(V,\mathbb{R}^{2}\right)$ for all $q\in \left [ 2,\infty  \right ] $, which are independent of $\lambda$. Namely, there exists a constants $C$ depending only on $q$ such that for any $\left(u,v\right) \in H_{\lambda}$,
$$\left \| \left ( u,v \right )  \right \| _{H}\leqslant C\left \| \left ( u,v \right ) \right \| _{H_{\lambda}},
\left \| \left ( u,v \right )  \right \| _{L^{q}\left(V,\mathbb{R}^{2}\right)}\leqslant C\left \| \left ( u,v \right ) \right \| _{H_{\lambda}}.$$
Moreover, the embedding of $H_{\lambda }$ into $L^{q}\left(V,\mathbb{R}^{2}\right)$ is compact for all $q\in \left [ 2,\infty  \right ] $. That is, for any bounded sequence $\left \{ \left ( u_{k},v_{k} \right )  \right \} \subset H_{\lambda}$, there exists $\left \{ \left ( u,v \right )  \right \} \subset H_{\lambda}$ such that, up to a subsequence,
\begin{equation}\label{convergesub}
\left\{\begin{array}{ll}
\left ( u_{k},v_{k} \right ) \rightharpoonup \left ( u,v \right )
&\textrm{in}  ~~H_{\lambda} \\
\left ( u_{k} ,v_{k}  \right )  \to  \left ( u ,v  \right )
&\forall x \in V\\
\left ( u_{k},v_{k} \right ) \to  \left ( u,v \right )
&\textrm{in}  ~~L^{q}\left(V,\mathbb{R}^{2}\right)
\end{array}\right.
\end{equation}
\end{lam}

\begin{proof}

Suppose that $\left ( u,v \right ) \in H_{\lambda}$, from $\left(A_{1}\right)$ one deduces immediately that $H_{\lambda}$ is embedded continuously into $H$.

For any vertex $x_{0}\in V$, we have
\begin{align*}
\left \| \left ( u,v \right )  \right \| _{H_{\lambda}}^{2} &=\int_{V}^{} \left ( \left | \nabla u \right | ^{2}+\left | \nabla v \right | ^{2}+\left ( \lambda a\left ( x \right ) +1 \right ) u^{2}+\left ( \lambda b\left ( x \right ) +1 \right ) v^{2} \right )d\mu \\
&\geqslant \int_{V}^{} \left ( \left | \nabla u \right | ^{2}+\left | \nabla v \right | ^{2}+u^{2}+v^{2} \right )d\mu \\
&\geqslant \int_{V}^{} u^{2} d\mu\\
&=\sum_{x\in V}^{} \mu \left ( x \right ) u^{2}\left ( x \right ) \\
&\geqslant \mu _{\min}u^{2}\left ( x_{0} \right )
\end{align*}
which gives $u\left ( x_{0} \right ) \leqslant\sqrt{\frac{1}{\mu_{\min}} }\left \| \left ( u,v \right )  \right \| _{H_{\lambda}} $ for any $x_{0} \in V$. Similarly, we have $v\left ( x_{0} \right ) \leqslant\sqrt{\frac{1}{\mu_{\min}} }\left \| \left ( u,v \right )  \right \| _{H_{\lambda}} $. Thus we have
\begin{displaymath}
\left \| \left ( u,v \right )  \right \| _{L^{\infty }}=\sup_{x\in V} \left | u\left ( x  \right ) \right | +\sup_{x\in V} \left | v\left ( x  \right ) \right | \leqslant2\sqrt{\frac{1}{\mu_{\min}} }\left \| \left ( u,v \right )  \right \| _{H_{\lambda}},
\end{displaymath}
which implies that $H_{\lambda} \hookrightarrow L^{\infty } \left (V,\mathbb{R}^{2}\right ) $ continuously and the embedding is independent of $\lambda$. Obviously, we also have $H_{\lambda} \hookrightarrow L^{2} \left (V,\mathbb{R}^{2}\right ) $ continuously. Then the interpolation gives the continuous embedding $H_{\lambda} \hookrightarrow L^{q } \left (V,\mathbb{R}^{2}\right ) $ for any $2\leqslant q \leqslant \infty $.

Since $H_{\lambda}$ is Hilbert space, for a bounded sequence $\left \{ \left ( u_{k},v_{k} \right )  \right \} $ in $H_{\lambda}$, there exists $\left ( u, v \right ) \in H_{\lambda} $ and a subsequence that we still call $\left \{ \left ( u_{k},v_{k} \right )  \right \} $, such that
\begin{align*}
\left ( u_{k},v_{k} \right ) \rightharpoonup \left ( u,v \right ) & ~~in~~H_{\lambda}\\
\left ( u_{k},v_{k} \right ) \rightharpoonup \left ( u,v \right ) & ~~in~~L^{2}\left ( V,\mathbb{R}^{2} \right ).
\end{align*}
Then we get that, for any $\left ( \xi ,\eta  \right ) \in L^{2}\left ( V,\mathbb{R}^{2} \right )$,
\begin{align}\label{weak}
& \lim_{k \to \infty}\int_{V}^{} \left(\left ( u_{k}-u   \right ) \xi +\left ( v_{k}-v \right ) \eta \right ) d \mu
\nonumber\\
= & \lim_{k \to \infty} \sum_{x\in V}^{} \left ( \mu \left ( x \right ) \left ( u_{k} -u \right ) \left ( x \right ) \xi \left ( x \right ) +\mu \left ( x \right ) \left ( v_{k} -v \right ) \left ( x \right ) \eta  \left ( x \right )  \right )
\nonumber\\
= & 0
\end{align}
Take any $x_{0} \in V $ and let
\begin{align}
\left ( \xi _{1},\eta_{1}  \right ) = &
\left \{
\begin{array}{ll}\label{L21}
\left ( 1,0 \right ) & x=x_{0}\\
\left (0,0 \right ) & x\neq x_{0}\\
\end{array}
 \right .
\\
\left ( \xi_{2} ,\eta_{2}  \right ) = &
\left \{
\begin{array}{ll}\label{L22}
\left ( 0,1 \right ) & x=x_{0}\\
\left (0,0 \right ) & x\neq x_{0}\\
\end{array}
\right .
\end{align}
Obviously, they both belong to $L^{2}\left ( V,\mathbb{R}^{2} \right )$.  By substituting (\ref{L21}) and (\ref{L22}) into (\ref{weak}), we get
\begin{align*}
\lim_{k \to \infty} \mu \left ( x_{0} \right ) \left ( u_{k}-u \right ) \left ( x_{0} \right ) =0\\
\lim_{k \to \infty} \mu \left ( x_{0} \right ) \left ( v_{k}-v \right ) \left ( x_{0} \right ) =0
\end{align*}
which implies that $\left ( u_{k} ,v_{k}  \right ) \to \left ( u ,v  \right )$ for all $x \in V$ while $k \to \infty$.

Next, we prove that $\left ( u_{k},v_{k} \right ) \to  \left ( u,v \right ) $ in $L^{q}\left ( V,\mathbb{R}^{2} \right )$. It is sufficient to prove that $\left ( u_{k},v_{k} \right ) \to  \left ( 0,0 \right ) $ strongly in $L^{q}\left ( V,\mathbb{R}^{2} \right )$ when $\left ( u_{k},v_{k} \right ) \rightharpoonup \left ( 0,0 \right ) $ weakly in $H_{\lambda}$.

Indeed, the boundedness of $\left \{ u_{k},v_{k} \right \} $ in $H_\lambda$ gives $\left \| \left ( u_{k},v_{k} \right )  \right \| _{H_{\lambda}}^{2} \leqslant C_{0}$ for some constant $C_{0}$. Since $a\left ( x \right )$ and $ b\left ( x \right ) $ satisfy $\left ( A_{2} \right ) $, for any $\epsilon > 0$, we can pick $R > 0$ such that $a\left ( x \right ) \geqslant \frac{2C_{0}}{\epsilon\lambda}$ and $b\left ( x \right ) \geqslant \frac{2C_{0}}{\epsilon\lambda} $ if $dist\left ( x,x_{0} \right ) > R$. Then we have
\begin{align}\label{convergeto0}
&\int_{dist\left ( x,x_{0} \right ) > R }^{} \left ( \left | u_{k} \right | ^{2}+\left | v_{k} \right | ^{2} \right )d\mu
\nonumber\\
\leqslant  &  \frac{\epsilon\lambda }{2C_{0}}\int_{dist\left ( x,x_{0} \right ) > R}^{}  \left ( a\left ( x \right ) \left | u_{k} \right | ^{2} + b\left ( x \right ) \left | v_{k} \right | ^{2} \right )d\mu
\nonumber\\
\leqslant  &  \frac{\epsilon }{2C_{0}}\left \| \left ( u_{k},v_{k} \right )  \right \|_{H_{\lambda}}^{2}
\nonumber\\
\leqslant  & \frac{\epsilon }{2}
\end{align}

On the other hand, since $\left \{ x \in V :  dist\left ( x,x_{0} \right ) \leqslant R \right \} $ is a finite set and $u_{k}\left ( x \right ) \to 0,v_{k}\left ( x \right ) \to 0$ for any $x \in V$ as $k \to \infty $, there exists $k_{0}$ such that $\int_{dist\left ( x,x_{0} \right ) \leqslant  R } \left ( \left | u_{k} \right | ^{2}+\left | v_{k} \right | ^{2} \right )d\mu
\leqslant \epsilon /2 $ when $k>k_{0}$. This together with (\ref{convergeto0}) gives that $\int_{V }^{} \left ( \left | u_{k} \right | ^{2}+\left | v_{k} \right | ^{2} \right )d\mu
\leqslant \epsilon$ when k is large enough. Therefore, $\lim_{k \to \infty} \left \| \left ( u_{k},v_{k} \right )  \right \| _{L^{2}}^{2} = 0$.

Since for $\left ( u_{k},v_{k} \right ) \in H_{\lambda}$ and any $x\in V$, we have
$\left \| \left ( u_{k},v_{k} \right )  \right \| _{L^{2}}^{2}\geqslant \mu_{min}u_{k}^{2} \left ( x \right ) ,
\left \| \left ( u_{k},v_{k} \right )  \right \| _{L^{2}}^{2}\geqslant \mu_{min}v_{k}^{2} \left ( x \right )$.
Hence,
\begin{align*}
\left \| \left ( u_{k},v_{k} \right )  \right \| _{L^{\infty }} & =\sup_{x\in V} \left | u_{k}\left ( x  \right ) \right | +\sup_{x\in V} \left | v_{k}\left ( x  \right ) \right |
\\
& \leqslant 2\sqrt{\frac{1}{\mu _{min}} } \left \| \left ( u_{k},v_{k} \right )  \right \| _{L^{2}}
\\
& \to 0 ~~as~~ k \to \infty
\end{align*}
Finally,for any $2<q< \infty$, one can get that
\begin{align*}
\left \| \left ( u_{k},v_{k} \right )  \right \| _{L^{q}}^{q} &
= \int_{V}^{}\left(\left | u_{k} \right |  ^{q}+\left | v_{k} \right |  ^{q}\right) d \mu
\\
& \leqslant \left [ \sup_{x \in V}\left | u_{k}\left ( x \right )  \right |  \right ] ^{q-2}\int_{V}^{}u_{k}^{2}\left ( x \right ) d\mu
+\left [ \sup_{x \in V}\left | v_{k}\left ( x \right )  \right |  \right ] ^{q-2}\int_{V}^{}v_{k}^{2}\left ( x \right ) d\mu
\\
& \leqslant \left \| \left ( u_{k},v_{k} \right )  \right \| _{L^{\infty }}^{q-2} \int_{V}^{} \left(u_{k}^{2}\left ( x \right ) +v_{k}^{2}\left ( x \right ) \right)d\mu
\\
& \to 0 ~~as~~ k \to \infty
\end{align*}
This completes the proof.
\end{proof}

For the function space $H_{\Omega}$, we also have a similar lemma.

\begin{lam}
Assume that $\Omega_{a}$ and $\Omega_{b}$ are bounded domains in $V$. Then $H_{\Omega}$ is continuously embedded into $L^{q_{1}}\left( \Omega_{a} \right)\times L^{q_{2}}\left( \Omega_{b} \right)$ for any $q_{1},q_{2} \in \left [1,\infty \right ] $. Namely, there exists a constant C depending only on $\Omega_{a}, \Omega_{b}$ and $q_{1},q_{2}$ such that for any $\left ( u,v \right ) \in H_{\Omega}$, $\left \| \left ( u,v \right )  \right \| _{L^{q_{1}}\left( \Omega_{a} \right)\times L^{q_{2}}\left( \Omega_{b} \right) } \leqslant C\left \| \left ( u,v \right )  \right \| _{H_{\Omega} } $. Moreover, for any bounded sequence $\left \{ \left ( u_{k},v_{k} \right )  \right \} \subset H_{\Omega}$, there exists $\left ( u,v \right ) \in H_{0}$ such that, up to a subsequence,

\begin{equation}
\left\{\begin{array}{ll}
\left ( u_{k},v_{k} \right ) \rightharpoonup \left ( u,v \right )
&\textrm{in}  ~~H_{\Omega}
\\
u_{k}\left ( x \right ) \to u\left ( x \right )
&\forall x \in \Omega _{a}
\\
v_{k}\left ( x \right ) \to v \left ( x \right )
&\forall x \in \Omega _{b}
\\
\left ( u_{k},v_{k} \right ) \to  \left ( u,v \right )
&\textrm{in}  ~~L^{q_{1}}\left( \Omega_{a} \right)\times L^{q_{2}}\left( \Omega_{b} \right)
\end{array}\right.
\end{equation}

\end{lam}
\begin{proof}
For $q_{1},q_{2} \in \left [2,\infty \right ] $, the proof is almost the same as that of Lemma \ref{embedding1}.

When $q_1=q_2=1$, since $\Omega_{a}$ and $\Omega_{b}$ are finite set, it is easy to prove that
\begin{displaymath}
 \left \| \left ( u,v \right )  \right \| _{L^{1}\left( \Omega_{a} \right)\times L^{1}\left( \Omega_{b} \right)}
=\sum_{x\in \Omega _{a}}\mu \left ( x \right )|u\left ( x \right )|
+\sum_{x\in \Omega _{b}}\mu \left ( x \right )|v\left ( x \right )|
\leqslant C \left \| \left ( u,v \right )  \right \|_{H_{\Omega}}
\end{displaymath}
where C is a constant depending on $\Omega_{a}$ and $\Omega_{b}$.
Moreover, noting that $u_{k}\left ( x \right ) \to u\left ( x \right ) $ in $\Omega_{a}$ and $v_{k}\left ( x \right ) \to v\left ( x \right ) $ in $\Omega_{b}$, for any given $\epsilon>0$, we have $\left | u_{k}\left ( x \right ) - u\left ( x \right ) \right | <\epsilon $ and $\left | v_{k}\left ( x \right ) - v\left ( x \right ) \right | <\epsilon $, where $k$ is large enough and independent of $x$. Thus,
\begin{align*}
\left \| \left ( u_{k}-u,v_{k}-v \right )  \right \| _{L^{1}\left( \Omega_{a} \right)\times L^{1}\left( \Omega_{b} \right)}
&=\sum_{x\in \Omega _{a}}\mu \left ( x \right )\left | u_{k}\left ( x \right )-u\left ( x \right )  \right |
+\sum_{x\in \Omega _{b}}\mu \left ( x \right )\left | v_{k}\left ( x \right )-v\left ( x \right )  \right |
\\
&\leqslant \epsilon \left ( \left | \Omega _{a} \right | +\left | \Omega _{b} \right | \right ) \to  0
\end{align*}
where $\left | \Omega _{a} \right |=\sum_{x\in \Omega _{a}}\mu \left ( x \right )$ and $\left | \Omega _{b} \right |=\sum_{x\in \Omega _{a}}\mu \left ( x \right )$ are bounded.

Therefore, we can verify that the embedding of $H_{\Omega}$ into $L^{q_{1}}\left( \Omega_{a} \right)\times L^{q_{2}}\left( \Omega_{b} \right)$ is compact for any $q_{1},q_{2} \in \left [1,\infty \right ] $.
\end{proof}

\section{Existence of solutions}

In this section, we first verify a compactness condition of the functional $J_\lambda$.

\begin{lam}\label{pscondition}
The functional $J_{\lambda}$ satisfies the $\left ( PS \right ) _{c}$ condition for every $c \in \mathbb{R}$. Namely, for any sequence $\left \{ \left ( u_{k},v_{k} \right )  \right \} \subset H_\lambda$ such that $J_\lambda\left ( u_{k},v_{k} \right )\to c $ and $J_\lambda '\left ( u_{k},v_{k} \right ) \to 0$, there is a convergent subsequence in $H_\lambda$.
\end{lam}
\begin{proof}
\quad
First, we claim that $\left \{ \left ( u_{k},v_{k} \right )  \right \} $ is bounded in $H_{\lambda}$.

In fact, assume that $\left \{ \left ( u_{k},v_{k} \right )  \right \}$ is a Palais-Smale sequence at level $c$. Obviously, $J_\lambda\left ( u_{k},v_{k} \right )\to c $ is equivalent to
\begin{equation}
\frac{1}{2}\left \| \left ( u_{k} ,v_{k}\right )  \right \| _{H_{\lambda }}^{2}
-\frac{1}{\alpha +\beta } \int _{V}\left | u_{k} \right | ^{\alpha }\left | v_{k} \right | ^{\beta }d\mu
=c+o_{k}\left (1\right ).
\label{ps1}
\end{equation}
Here and in the sequel, $o_{k}\left(1\right)\to 0$ as $k\to +\infty$.
Since $J_\lambda '\left ( u_{k},v_{k} \right ) \to 0$, we get
\begin{align*}
\left \langle J_{\lambda }'\left ( u_{k},v_{k} \right ) ,\left ( u_k ,v_k  \right )  \right \rangle
&=\left \langle \left ( u_{k},v_{k} \right ) ,\left ( u_k ,v_k  \right )  \right \rangle_{H_\lambda}
-\int _{V}\left ( \frac{\alpha }{\alpha +\beta }\left | u_{k} \right | ^{\alpha -2} u_{k} \left | v_{k} \right | ^{\beta } u_k
+\frac{\beta  }{\alpha +\beta }\left | u_{k} \right | ^{\alpha}\left | v_{k} \right | ^{\beta -2}v_{k} v_k
\right ) d\mu
\notag\\
&=o_{k}\left(1\right)\left \| \left ( u_k ,v_k  \right )  \right \| _{H_{\lambda }},
\end{align*}
which is equivalent to
\begin{equation}
\left \| \left ( u_{k},v_{k}  \right )  \right \| _{H_{\lambda }}^{2}
=\int _{V}\left | u_{k} \right |^{\alpha }\left | v_{k} \right |^{\beta  }d\mu
+o_{k}\left (1\right ) \left \| \left ( u_{k},v_{k}  \right )  \right \| _{H_{\lambda }}
\label{ps3}
\end{equation}
Then combining \eqref{ps1} and \eqref{ps3}, we conclude that
\begin{equation}\label{bounded}
 \left (\frac{1}{2}-\frac{1}{\alpha +\beta }  \right ) \left \| \left ( u_{k},v_{k}  \right )  \right \| _{H_{\lambda }}^{2}
=c+ o_{k}\left (1\right )
+o_{k}\left (1\right ) \left \| \left ( u_{k},v_{k}  \right )  \right \| _{H_{\lambda }}
\end{equation}
Since $\alpha >1$ and $\beta >1$, we obtain that $\left \{ \left ( u_{k},v_{k} \right )  \right \} $ is bounded. By Lemma \ref{embedding1}, there exists $\left ( u,v \right ) \in H_\lambda$ and a subsequence, still denoted by $\left \{ \left ( u_{k},v_{k} \right )  \right \} $, such that satisfying (\ref{convergesub}).

We now show that the convergence of $\left ( u_{k},v_{k} \right )$ to $ \left ( u,v \right )$ is strong in $H_\lambda$ which implies the $(PS)_c$ condition. First, we have
\begin{align*}
& \int_{V}^{} \left | u_{k} \right |^{\alpha -2} u_{k}\left | v_{k} \right | ^{\beta} \left ( u_{k}-u \right ) d\mu
\notag\\
\leqslant & \left [ \sup_{x \in V} \left | u_{k} \right | \right ] ^{\alpha -1}\left ( \int_{V}^{} \left | v_{k} \right | ^{2\beta } d\mu \right ) ^{\frac{1}{2} }
\left ( \int_{V}^{} \left | u_{k}-u \right | ^{2} d\mu \right ) ^{\frac{1}{2} }
\notag\\
= & \left \| u_{k} \right \| _{L^{\infty }}^{\alpha -1}\left \| v_{k} \right \| _{L^{2\beta}}^{\beta} \left \| u_{k}-u \right \| _{L^{2}}
\notag\\
\to &~0~~as ~~k \to \infty
\end{align*}
Similarly, we also have $\int_{V}^{} \left | u_{k} \right | ^{\alpha } \left | v_{k} \right | ^{\beta -2}v_{k}\left ( v_{k}-v \right ) d \mu \to 0$ as $k \to \infty$.

$J_\lambda '\left ( u_{k},v_{k} \right ) \to 0$ leads to $\left \langle J_{\lambda }'\left ( u_{k},v_{k} \right ) ,\left(u_{k}-u,v_{k}-v\right)  \right \rangle \to 0$, which is equivalent to
\begin{align}
&\left \langle \left(u_{k},v_{k}\right) ,\left(u_{k}-u,v_{k}-v\right) \right \rangle_{H_\lambda}
\notag \\
=&\int _{V}\left ( \frac{\alpha }{\alpha +\beta }\left | u_{k} \right | ^{\alpha -2} u_{k} \left | v_{k} \right | ^{\beta } \left(u_{k}-u\right)
+\frac{\beta  }{\alpha +\beta }\left | u_{k} \right | ^{\alpha}\left | v_{k} \right | ^{\beta -2}v_{k} \left(v_{k}-v\right)
\right ) d\mu
\notag \\
&+o_{k}\left(1\right)\left \| \left ( u_{k}-u,v_{k}-v \right )  \right \| _{H_{\lambda }}
\notag\\
\to &~0~~as ~~k \to \infty
\label{ps5}
\end{align}
Moreover, as $\left ( u_{k},v_{k} \right ) \rightharpoonup \left ( u,v \right ) $ weakly in $H_{\lambda}$, there holds
\begin{equation}
\left \langle \left(u,v\right) ,\left(u_{k}-u,v_{k}-v\right) \right \rangle_{H_\lambda} \to ~0~~as ~~k \to \infty
\label{ps6}
\end{equation}

Combining \eqref{ps5} and \eqref{ps6}, we get that $\left ( u_{k},v_{k} \right ) \to \left ( u,v \right ) $ strongly in $H_{\lambda}$ and the lemma is proved.
\end{proof}

To prove the existence results, we need to check the geometric conditions of mountain pass theorem, which are presented in the following lemma.

\begin{lam}\label{mountain}
The function $J_{\lambda}$ satisfies the mountain pass geometry. Namely,

$\left ( \rmnum{1} \right ) $ there exist positive constants $r$ and $\rho$, such that $J_{\lambda}\left ( u,v \right ) > r$ for $\left \| \left ( u,v \right )  \right \|_{H_\lambda} =\rho $

$\left ( \rmnum{2} \right ) $ there exists $\left( u_0, v_0 \right) \in H_\lambda \backslash \left \{ \left ( 0,0 \right )  \right \} $ such that $\left \| \left( u_0, v_0 \right) \right \|_{H_\lambda} > \rho $ and $J_{\lambda} \left ( u_0, v_0 \right ) <0$.
\end{lam}
\begin{proof}
\quad
 First, for any $\left ( u,v \right ) \in H_{\lambda}$, we have
\begin{align*}
J_{\lambda} \left ( u,v \right ) & = \frac{1}{2} \left \| \left ( u,v \right )  \right \| _{H_\lambda}^{2} -\frac{1}{\alpha+\beta} \int_{V}^{} \left | u \right | ^{\alpha } \left | v \right | ^{\beta } d\mu
\\
& \geqslant \frac{1}{2} \left \| \left ( u,v \right )  \right \| _{H_\lambda}^{2} -\frac{1}{\alpha+\beta} \int_{V}^{}\left( \left | u \right | ^{\alpha +\beta } + \left | v \right | ^{\alpha +\beta }\right) d\mu
\\
& =\frac{1}{2} \left \| \left ( u,v \right )  \right \| _{H_\lambda}^{2}-\frac{1}{\alpha+\beta}\left \| \left ( u,v \right )  \right \| _{L^{\alpha +\beta }}^{\alpha +\beta }
\\
& \geqslant \frac{1}{2} \left \| \left ( u,v \right )  \right \| _{H_\lambda}^{2}-C\left \| \left ( u,v \right )  \right \| _{H_\lambda}^{\alpha +\beta }
\end{align*}
Since $\alpha +\beta >2$, we can choose an enough small constant $\rho$, such that for any $\left ( u,v \right )$ with its $H_\lambda$ norm equals to $\rho$, there holds $J_{\lambda} \left ( u,v \right ) \geqslant\frac{1}{2} \rho ^{2}-C \rho ^{\alpha +\beta }>0$. Part $\left ( \rmnum{1} \right ) $ is proved.

On the other hand, for any $\left ( u,v \right )  \in H_{\lambda} \backslash \left \{ \left ( 0,0 \right )  \right \}$ and $t>0$ we have
\begin{displaymath}
J_{\lambda}\left ( \left ( tu,tv \right )  \right ) =\frac{t^{2}}{2} \left \| \left ( u,v \right )  \right \| _{H_\lambda}^{2}-\frac{t^{\alpha +\beta }}{\alpha+\beta} \int_{V}^{} \left | u \right | ^{\alpha } \left | v \right | ^{\beta }d \mu
\end{displaymath}
Therefore $J_{\lambda}\left ( \left ( tu,tv \right )  \right ) \to -\infty $ as $t \to + \infty $,  and we can certainly find $\left( u_0, v_0 \right) \in H_\lambda $ such that $\left \| \left( u_0, v_0 \right) \right \|_{H_\lambda} \geqslant \rho $ and $J_{\lambda} \left( u_0, v_0 \right) <0$.
\end{proof}

\begin{lam}\label{mpcritical}
The functional $J_{\lambda}$ has a nontrivial critical point.
\end{lam}
\begin{proof}
For each $\lambda>0$, by Lemma \ref{mountain}, $J_{\lambda}$ satisfies all the hypotheses of the mountain-pass geometry. Let $\left( u_0, v_0 \right)$ and $r$ be given by the lemma. We may define the mountain pass level $c_\lambda$ of $J_{\lambda}$ as
\begin{displaymath}
 c_{\lambda }:=\inf _{\gamma \in \Gamma }\sup _{t\in \left [ 0,1 \right ] }J_{\lambda }
\left ( \gamma \left ( t \right )  \right ) \geqslant r >0
\end{displaymath}
where
\begin{displaymath}
 \Gamma :=\left \{ \gamma \in C\left ( \left [ 0,1 \right ] ,H_{\lambda } \right ) :\gamma \left ( 0 \right ) =\left ( 0,0 \right ) ,\gamma \left ( 1 \right ) =\left( u_0, v_0 \right) \right \}
\end{displaymath}
Then the mountain pass theorem gives that there exists $\left \{ \left ( u_{k},v_{k} \right )  \right \} \subset H_{\lambda }$ such that $J_{\lambda }\left ( u_{k},v_{k} \right ) \to c_{\lambda }$ and $J_{\lambda }'\left ( u_{k},v_{k} \right ) \to 0$. Moreover, by Lemma \ref{pscondition}, $J_{\lambda}$ satisfies the Palais-Smale condition. Namely, $\left \{ \left ( u_{k},v_{k} \right )  \right \}$ is bounded in $H_{\lambda}$ and, up to a subsequence, we may assume that $ \left ( u_{k},v_{k} \right )$ converge to some $\left ( u_{\lambda },v_{\lambda } \right )$ strongly in $H_{\lambda}$.

Then it is easy to conclude that, for any $\left(\xi,\eta \right) \in C_{c}\left(V\right)\times C_{c}\left(V\right)$, we have
\begin{displaymath}
 \left \langle J_{\lambda }'\left ( u_{\lambda },v_{\lambda } \right ) ,\left ( \xi ,\eta \right )  \right \rangle
=\lim_{k\to +\infty } \left \langle J_{\lambda }'\left ( u_{k},v_{k } \right ) ,\left ( \xi ,\eta \right )  \right \rangle
=0,
\end{displaymath}
which implies that $\left ( u_{\lambda },v_{\lambda } \right )$ is a weak solution of the system (\ref{system}). In addition, $c_{\lambda}$ is a critical level of the functional $J_{\lambda}$, which is achieved at $\left(u_{\lambda},v_{\lambda}\right)$. Since $c_{\lambda}=J_{\lambda}\left(u_{\lambda},v_{\lambda}\right)>0$, we conclude that $\left(u_{\lambda},v_{\lambda}\right)$ is a nontrivial critical point of $J_{\lambda}$.
\end{proof}

In fact, the critical point $\left(u_{\lambda},v_{\lambda}\right)$ given by the above lemma is a ground state solution. To prove this fact, we introduce the Nehari manifold associated to system (\ref{system}), which is defined by
\begin{displaymath}
 \mathcal{N} _{\lambda }:=\left \{ \left ( u,v \right ) \in H_{\lambda } \setminus   \left \{ \left ( 0,0 \right )  \right \} :
\left \langle J_{\lambda }'\left ( u,v \right ) ,\left ( u,v \right )  \right \rangle =0 \right \}
\end{displaymath}
The least energy level of the functional among the Nehari manifold is
\begin{displaymath}
 c_{\mathcal{N}_{\lambda } }=\inf _{\left ( u,v \right )\in \mathcal{N}_{\lambda } }J_{\lambda }\left ( u,v \right )
\end{displaymath}
The following lemma tells us that the least energy level and the mountain pass level are the same.

\begin{lam}
The two levels $c_{\lambda}$ and $c_{\mathcal{N}_{\lambda } }$ satisfy $c_{\lambda}=c_{\mathcal{N}_{\lambda } }$.
\end{lam}
\begin{proof}
For the nontrivial critical point $\left(u_{\lambda},v_{\lambda}\right)$ of $J_{\lambda}$ given by Lemma \ref{mpcritical}, we have that $\left(u_{\lambda},v_{\lambda}\right) \neq \left(0,0\right)$ and $\left \langle J_{\lambda }'\left ( u_{\lambda },v_{\lambda } \right ) ,\left ( u_{\lambda },v_{\lambda } \right ) \right \rangle =0$. It is obviously that $\left ( u_{\lambda },v_{\lambda } \right ) \in \mathcal{N}_{\lambda } $, which implies that $c_{\mathcal{N}_{\lambda } }\leqslant c_{\lambda}=J_{\lambda }\left ( u_{\lambda },v_{\lambda } \right ) $.

Now it is sufficient to show that $c_{\lambda}\leqslant c_{\mathcal{N}_{\lambda } }$. Fix any $\left(u,v\right)\in \mathcal{N}_{\lambda }$, in view of the proof of Lemma \ref{mountain}, there exists some $t_{0}>0$ large enough such that $J_{\lambda}\left(t_{0}u,t_{0}v\right)<0$. Therefore, by taking $\left(u_{0},v_{0}\right)=\left(t_{0}u,t_{0}v\right)$, we can define $\gamma_0 :\left [ 0,1 \right ] \to H_{\lambda }$ as $\gamma_0 \left(t\right)=\left(tu_{0},tv_{0}\right)$, which satisfies $\gamma_0 \in \Gamma$. Recalling the definition of $c_{\lambda}$, we can conclude that
\begin{align*}
c_{\lambda }\leqslant \sup_{t\in \left [ 0,1 \right ] }J_{\lambda }\left ( \gamma_0 \left ( t \right )  \right )
=\sup_{t\in \left [ 0,1 \right ] }J_{\lambda }\left ( tt_{0}u,tt_{0}v \right )
\leqslant  \sup_{t\geqslant0 }J_{\lambda }\left ( tu,tv \right )
\end{align*}
Differentiating $J_{\lambda }\left ( tu,tv \right )$ with respect to $t$ and noticing that $\left \langle J_{\lambda }'\left ( u,v \right ) ,\left ( u,v \right ) \right \rangle =0$, we get
\begin{displaymath}
 \frac{dJ_{\lambda }\left ( tu,tv \right ) }{dt} =\left ( t-t^{\alpha +\beta -1} \right )
\left \| \left ( u,v \right )  \right \|  _{H_{\lambda }}^{2}.
\end{displaymath}
It is not hard to check that $J_{\lambda }\left ( tu,tv \right )$ arrives its maximum value when $t=1$, which implies that $\sup_{t\geqslant 0 }J_{\lambda }\left ( tu,tv \right )=J_{\lambda }\left ( u,v \right )$. According to the arbitrariness of $\left(u,v\right)\in \mathcal{N}_{\lambda } $, we can obtain that $c_{\lambda }\leqslant \inf _{\left ( u,v \right )\in \mathcal{N}_{\lambda } }J_{\lambda }\left ( u,v \right ) =c_{\mathcal{N}_{\lambda } }$.
\end{proof}

Similarly, the Nehari manifold associated to system (\ref{dirichlet}) is defined by
\begin{displaymath}
 \mathcal{N} _{\Omega }:=\left \{ \left ( u,v \right ) \in H_{\Omega } \setminus   \left \{ \left ( 0,0 \right )  \right \} :
\left \langle J_{\Omega }'\left ( u,v \right ) ,\left ( u,v \right )  \right \rangle =0 \right \}
\end{displaymath}
The least energy level of the functional $J_{\Omega }$ among the Nehari manifold $\mathcal{N} _{\Omega }$ is
\begin{displaymath}
 c_{\mathcal{N}_{\Omega } }=\inf _{\left ( u,v \right )\in \mathcal{N}_{\Omega } }J_{\Omega }\left ( u,v \right ).
\end{displaymath}
Now we are ready to prove Theorem \ref{ssolution} and \ref{dsolution}.

\begin{proof}
From the above arguments, we can conclude that for any fixed $\lambda>0$, $J_{\lambda}$ has a nontrivial critical point $\left(u_{\lambda},v_{\lambda}\right)$ at level $c_{\lambda}$. Moreover, $c_{\lambda}$ is a ground state of system (\ref{system}). This completes the proof of Theorem \ref{ssolution}.

Since finite graphs are a specialization of locally finite graphs, the proofs of the previous results can be easily applied to the system (\ref{dirichlet}) with minor modifications. The corresponding results are still valid and we can get a ground state solution $(u_\Omega,v_\Omega)\in \mathcal{N} _{\Omega }$ of the system (\ref{dirichlet}) at level $c_{\mathcal{N}_{\Omega } }$. We omit the details of the proofs here.
\end{proof}

\begin{Rem}\label{level}
From (\ref{bounded}) and the weak lower semi-continuity of the norm $\|\cdot\|_{H_\lambda}$, we can deduce that
$$
\|(u_\lambda,v_\lambda)\|^2_{H_\lambda}\leqslant \frac{2(\alpha+\beta)c_{\mathcal{N}_{\lambda } }}{\alpha+\beta-2}.
$$
Therefore, $\|(u_\lambda,v_\lambda)\|_{H_\lambda}^2$ is uniformly bounded by the constant $\frac{2(\alpha+\beta)c_{\mathcal{N}_{\Omega } }}{\alpha+\beta-2}$, which is independent of $\lambda$.
\end{Rem}

\section{The asymptotic behavior}

We devote this section to the proof of Theorem \ref{convergence}. Suppose that $\left(u_{\lambda_{n}},v_{\lambda_{n}}\right)$ is the nontrivial solution of the system (\ref{system}) obtained according to Theorem \ref{ssolution}, where $\lambda_{n}$ tends to $\infty$ as $n\to \infty$. For simplicity, from now on we use $u_{n}$ and $v_{n}$ to denote $u_{\lambda_{n}}$ and $v_{\lambda_{n}}$ respectively. Since $\left\|\left(u_{n},v_{n}\right)\right\|_H\leqslant \left\|\left(u_{n},v_{n}\right)\right\|_{H_\lambda}$, Remark \ref{level} tells us that $\left(u_{n},v_{n}\right)$ is uniformly bounded both in $H$ and $H_\lambda$. Up to a subsequence, there exists some $\left(u,v\right)\in H$ such that
\begin{displaymath}
\left\{\begin{array}{ll}
\left ( u_{n},v_{n} \right ) \rightharpoonup \left ( u,v \right )
&\textrm{in}  ~~H\\
\left ( u_{n}\left ( x \right ) ,v_{n}\left ( x \right )  \right )  \to  \left ( u\left ( x \right ) ,v\left ( x \right )  \right )
&\forall x \in V\\
\left ( u_{n},v_{n} \right ) \to  \left ( u,v \right )
&\textrm{in}  ~~L^{q}\left(V,\mathbb{R}^{2}\right)
\end{array}\right.
\end{displaymath}

Firstly, we claim that $u\equiv 0$ in $\Omega_{a}^{c}$. If not, we can find some vertex $x_{0}\in \Omega_{a}^{c}$ such that $u\left(x_{0}\right) \neq 0$ and get
\begin{displaymath}
\frac{1}{\lambda _{n}}\left \|(u_{n},v_n) \right \|_{H_{\lambda _{n}}}^{2}
\geqslant
\frac{1}{\lambda _{n}}\int _{V}\lambda _{n}a \left (x\right )u_{n}^{2}d\mu
\geqslant
u_{n}^{2}\left (x_{0}\right )a\left (x_{0}\right ) \mu \left (x_{0}\right )
>0.
\end{displaymath}
Recall that $(u_{n},v_{n})$ is bounded by a constant independent of $\lambda$, so we have that $\frac{1}{\lambda _{n}}\left \|(u_{n},v_n) \right \|_{H_{\lambda _{n}}}^{2} \to 0$. On the other hand, we have $u_{n}^{2}\left (x_{0}\right )a\left (x_{0}\right ) \mu \left (x_{0}\right )\to u^{2}\left (x_{0}\right )a\left (x_{0}\right ) \mu \left (x_{0}\right )>0$ as $n\to \infty$, which is a contradiction. Then the claim is proved and we can conclude that $u\in W_{0}^{1,2}\left(\Omega_{a}\right)$. Analogously, we also have $v\in W_{0}^{1,2}\left(\Omega_{b}\right)$.

Next, we prove that $(u,v)$ is a solution of the system (\ref{dirichlet}). For any given $\varphi \in C_{c}\left(\Omega_{a}\right)$, using $\left(\varphi , 0 \right)$ as a test function, we have $\left \langle J'_\lambda\left ( u_{n},v_{n} \right ) ,\left ( \varphi  ,0 \right )  \right \rangle =0$. Namely,
\begin{displaymath}
\int _{V}\left(\nabla u_{n}\nabla \varphi +\left (\lambda a\left (x\right )+1\right )u_{n}\varphi\right) d\mu
=\int _{V}\frac{\alpha }{\alpha +\beta }\left |u_{n}\right |^{\alpha -2}u_{n}\left |v_{n}\right |^{\beta } \varphi d\mu
\end{displaymath}
Since $\varphi=0$ on $\Omega_{a}^{c}$ and $a\left(x\right)=0$ on $\Omega_{a}$, we obtain
\begin{equation}\label{testfunc}
\int _{\partial \Omega _{a}\cup  \Omega _{a}}\nabla u_{n}\nabla \varphi d\mu +\int _{\Omega _{a}}u_{n}\varphi d\mu
=\int _{\Omega _{a}}\frac{\alpha }{\alpha +\beta }\left |u_{n}\right |^{\alpha -2}u_{n}\left |v_{n}\right |^{\beta } \varphi d\mu
\end{equation}
The weakly convergence of $u_{n}$ to $u$ in $H$ gives $\left \langle u_{n},\varphi  \right \rangle_H \to \left \langle u,\varphi  \right \rangle_H $. And the pointwise convergence of $\left ( u_{n}\left ( x \right ) ,v_{n}\left ( x \right )  \right )$ to $\left ( u\left ( x \right ) ,v\left ( x \right )  \right )$ in $V$ gives
\begin{displaymath}
\int _{\Omega _{a}}\left |u_{n}\right |^{\alpha -2}u_{n}\left |v_{n}\right |^{\beta }\varphi d\mu
\to \int _{\Omega _{a}}\left |u\right |^{\alpha -2}u\left |v\right |^{\beta }\varphi d\mu
\end{displaymath}
Therefore,  as $n\to \infty$, the equation (\ref{testfunc}) turns to
\begin{displaymath}
\int _{\partial \Omega _{a}\cup  \Omega _{a}}\nabla u\nabla \varphi d\mu +\int _{\Omega _{a}}u\varphi d\mu
=\int _{\Omega _{a}}\frac{\alpha}{\alpha+\beta}\left |u\right |^{\alpha -2}u\left |v\right |^{\beta }\varphi d\mu
\end{displaymath}
Since $\nabla \varphi=0$ on $\left (\overline{\Omega }_{a}\right )^{c} $, the above equation is equivalent to
\begin{equation}\label{func1}
\int _{\overline{\Omega }_{a}\cup \overline{\Omega }_{b}}\nabla u_{n}\nabla \varphi d\mu
+\int _{\Omega _{a}\cup \Omega _{b}}u_{n}\varphi d\mu
=\int _{\Omega _{a}\cup \Omega _{b}}\frac{\alpha}{\alpha+\beta}\left |u\right |^{\alpha -2}u\left |v\right |^{\beta }\varphi d\mu
\end{equation}
Analogously, for any given $\psi \in C_{c}\left(\Omega_{b}\right)$ and using $\left(0,\psi \right)$ as a test function, we have
\begin{equation}\label{func2}
\int _{\overline{\Omega }_{a}\cup \overline{\Omega }_{b}}\nabla v_{n}\nabla \psi d\mu
+\int _{\Omega _{a}\cup \Omega _{b}}v_{n}\psi d\mu
=\int _{\Omega _{a}\cup \Omega _{b}}\frac{\beta}{\alpha+\beta}\left |v\right |^{\alpha }\left |v\right |^{\beta -2}v \psi d\mu
\end{equation}
Combining (\ref{func1}) and (\ref{func2}),  we conclude that for any $\left(\varphi,\psi \right) \in C_{c}\left(\Omega_{a}\right) \times C_{c}\left(\Omega_{b}\right)$, $\left \langle J'_\Omega\left ( u,v \right ) ,\left ( \varphi ,\psi  \right )  \right \rangle =0$, which tells us that $(u,v)$ is a critical point of the functional $J_\Omega$.

We now prove that $(u, v)$ is a nontrivial solution of the system (\ref{dirichlet}). Since $\left(u_{n},v_{n}\right)$ is a nontrivial critical point of $J_{\lambda}$, we have
\begin{align*}
0 & =\left \langle J'_{\lambda }\left ( u_{n},v_{n} \right ) ,\left ( u_{n},v_{n} \right )  \right \rangle
\\
&=\left \| \left ( u_{n},v_{n} \right )  \right \| _{H_{\lambda }}^{2}
-\int _{V}\left | u_{n} \right | ^{\alpha }\left | v_{n} \right | ^{\beta }d\mu
\\
&\geqslant \left \| \left ( u_{n},v_{n} \right )  \right \| _{H_{\lambda }}^{2}
-\int _{V}\left(\left | u_{n} \right | ^{\alpha +\beta }+\left | v_{n} \right | ^{\alpha +\beta }\right)d\mu
\\
&= \left \| \left ( u_{n},v_{n} \right )  \right \| _{H_{\lambda }}^{2}
-\left \| \left ( u_{n},v_{n} \right )  \right \| _{L^{\alpha +\beta }}^{\alpha +\beta }
\\
&\geqslant \left \| \left ( u_{n},v_{n} \right )  \right \| _{H_{\lambda }}^{2}
-C\left \| \left ( u_{n},v_{n} \right )  \right \| _{H_{\lambda }}^{\alpha +\beta },
\end{align*}
where $C$ in the last inequality is the embedding constant in Lemma \ref{embedding1}. Then for the constant $\sigma:=\left ( \frac{1}{C}  \right ) ^{\frac{1}{\alpha +\beta -2} }$, we have that $\left \| \left ( u_{n}, v_{n} \right )  \right \| _{H_{\lambda }}\geqslant \sigma>0 $.
Since $(u_n,v_n)\in \mathcal{N}_{\lambda }$, we have
\begin{equation*}
\sigma^2\leqslant \left \| \left ( u_{n}, v_{n} \right )  \right \|^2 _{H_{\lambda }}
=\int_{V}|u_n|^\alpha |v_n|^\beta d\mu
\leqslant \int_V \left(|u_n|^{\alpha+\beta}+|v_n|^{\alpha+\beta}\right)d\mu
\end{equation*}
Then from $\left( u_{n},v_{n} \right ) \to  \left ( u,v \right )$ in $L^{q}\left(V,\mathbb{R}^{2}\right)$, we can conclude that $\left(u,v\right) \neq \left(0,0\right)$ and consequently we have $\left(u,v\right)\in \mathcal{N}_{\Omega}$.

To prove Theorem \ref{convergence}, we also need to verify that $\left(u,v\right)$ achieves the infimum of $J_{\Omega}$ in $\mathcal{N}_{\Omega}$. Noticing that $\mathcal{N} _{\Omega }$ is a subset of $\mathcal{N} _{\lambda }$, we have that $c_{\mathcal{N}_{\lambda } }\leqslant c_{\mathcal{N}_{\Omega } }$ for any $\lambda>0$. Therefore, for any $\lambda_n$, we have
\begin{align*}
c_{\Omega }\geqslant c_{\lambda_n }=J_{\lambda_n }\left ( u_{n},v_{n} \right )
& = J_{\lambda_n }\left ( u_{n},v_{n} \right ) -\frac{1}{2}\left \langle J'_{\lambda_n } \left ( u_{n},v_{n} \right ),\left ( u_{n},v_{n} \right )\right \rangle
\\
&=\left ( \frac{1}{2}-\frac{1}{\alpha +\beta }\right ) \int _{V} \left | u_{n} \right | ^{\alpha }\left | v_{n} \right | ^{\beta } d\mu,
\end{align*}
which implies that
\begin{align}\label{fatou}
c_{\Omega }\geqslant\limsup _{\lambda _{n}\to \infty } c_{\lambda_{n} }
&=\limsup _{\lambda _{n}\to \infty }\left ( \frac{1}{2}-\frac{1}{\alpha +\beta } \right )\int _{V} \left | u_{n} \right | ^{\alpha }\left | v_{n} \right | ^{\beta } d\mu\nonumber
\\
&\geqslant\left ( \frac{1}{2}-\frac{1}{\alpha +\beta } \right )\int _{\Omega _{a}\cup \Omega _{b}} \left | u \right | ^{\alpha }\left | v \right | ^{\beta } d\mu\nonumber
\\
&=J_\Omega\left ( u,v \right ) -\frac{1}{2} \left \langle J'_\Omega\left ( u,v \right ) ,\left ( u,v \right ) \right \rangle\nonumber
\\
&=J_\Omega\left ( u,v \right ) \geqslant c_{\Omega }
\end{align}
Hence, we have $J_\Omega\left ( u,v \right ) = c_{\Omega }$.

Finally, we show the strong convergence of $(u_n,v_n)$ to $(u,v)$ in $H$ to finish the proof of the theorem. First, noticing that $H_\lambda$ is a Hilbert space, the weak convergence of $(u_n,v_n)$ to $(u,v)$ in $H_\lambda$ gives
\begin{align}\label{strong1}
\|(u_n,v_n)-(u,v)\|^2_{H_\lambda}=\|(u_n,v_n)\|^2_{H_\lambda}-\left \langle (u,v), (u_n,v_n)\right \rangle_{H_\lambda}=\|(u_n,v_n)\|^2_{H_\lambda}-\|(u,v)\|^2_{H_\lambda}+o_n(1).
\end{align}
Since $(u_n,v_n)\in \mathcal{N}_{\lambda }$ and $(u,v)\in \mathcal{N}_{\Omega }$, we have
\begin{align}\label{strong2}
\|(u_n,v_n)\|^2_{H_\lambda}=\int_V |u_n|^\alpha |v_n|^\beta d\mu
\end{align}
and
\begin{align}\label{strong3}
\|(u,v)\|^2_{H_\lambda}=\|(u,v)\|^2_{H_\Omega}=\int_{\Omega_a\cup \Omega_b}|u|^\alpha |v|^\beta d\mu.
\end{align}
Combining (\ref{strong1}), (\ref{strong2}) and (\ref{strong3}), we get
\begin{align}\label{strong4}
\|(u_n,v_n)-(u,v)\|^2_{H}\leqslant \|(u_n,v_n)-(u,v)\|^2_{H_\lambda}
=\int_V |u_n|^\alpha |v_n|^\beta d\mu-\int_{\Omega_a\cup \Omega_b}|u|^\alpha |v|^\beta d\mu+o_n(1),
\end{align}
Since (\ref{fatou}) implies that
\begin{align*}
\int_V |u_n|^\alpha |v_n|^\beta d\mu-\int_{\Omega_a\cup \Omega_b}|u|^\alpha |v|^\beta d\mu=o_n(1),
\end{align*}
consequently, (\ref{strong4}) gives
\begin{align*}
\|(u_n,v_n)-(u,v)\|^2_{H}=o_n(1)
\end{align*}
and the proof is finished.

\section{Numerical results}

In order to better illustrate our conclusions, we design a finite connected graph $G_{22}=\left(V,E\right)$ who has 22 vertices. For simplicity, the positive measure on a vertex in $V$ and the symmetric weight of an edge in $E$ are all set to equal $1$. The structure of $G_{22}$ is shown in Figure 1. Obviously, the graph $G_{22}$ satisfies all the assumptions in the previous theorems. Furthermore, for the systems (\ref{system}) and (\ref{dirichlet}), we let $\alpha = \beta = 2$ and the potential functions are defined as follows.
\begin{displaymath}
 a\left (x_{i}\right )=
 \left \{
\begin{array}{ll}
0,& i = 1,2,3,4,5,6,7,8,9
\\
1, & i \neq 1,2,3,4,5,6,7,8,9
\end{array}
\right .
\end{displaymath}
\begin{displaymath}
 b\left (x_{i}\right )=
 \left \{
\begin{array}{ll}
0,& i = 1,2,3,4,5,6,10,11,12
\\
1, & i \neq 1,2,3,4,5,6,10,11,12
\end{array}
\right .
\end{displaymath}
With the definitions of the potential functions, we have
\begin{displaymath}
\begin{array}{cc}
 \Omega_{a}=\left\{x_{1},x_{2},x_{3},x_{4},x_{5},x_{6},x_{7},x_{8},x_{9}\right\}      \\
 \Omega_{b}=\left\{x_{1},x_{2},x_{3},x_{4},x_{5},x_{6},x_{10},x_{11},x_{12}\right\}  \\
 \Omega_{a}\cap \Omega_{b} = \left\{x_{1},x_{2},x_{3},x_{4},x_{5},x_{6}\right\}
\end{array}
\end{displaymath}
Their boundaries are
\begin{displaymath}
\begin{array}{cc}
\partial \Omega_{a}=\left\{x_{10},x_{11},x_{12},x_{13},x_{14},x_{16},x_{17},x_{18} ,x_{22}\right\}      \\
\partial \Omega_{b}=\left\{x_{7},x_{8},x_{9},x_{13},x_{17},x_{18},x_{19},x_{21} ,x_{22}\right\}
\end{array}
\end{displaymath}

\begin{figure}[htbp]
\centering
\includegraphics[width=3.6in,height=1.8in]{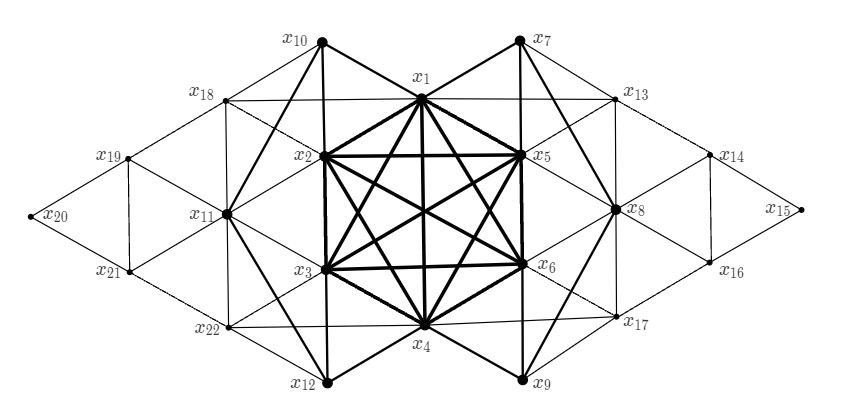}
\caption{The graph $G_{22}$}
\label{fig:graph}
\end{figure}

We compute the numerical solutions of the systems (\ref{system}) and (\ref{dirichlet}) by MATLAB. The numerical solution of the limit system (\ref{dirichlet}) is shown in Table 1 and the values not listed in the table are equal to zero. Henceforth, we use $u_{i}$ and $v_{i}$ to denote the value of the functions $u(x)$ and $v(x)$ at vertex $x_i$, where $i = 1,2,\cdots , 22$.

\begin{table}[!hbp]
    \centering
    \caption{The numerical solution of system (2)}
    \begin{tabular}{ccccccccc}
    \hline
    $u_{1}$&$u_{2}$&$u_{3}$&$u_{4}$&$u_{5}$&$u_{6}$&$u_{7}$&$u_{8}$&$u_{9}$ \\
    \hline
    3.5308   &  2.0210   &  2.0210   &  3.5308   &  2.1900   &  2.1900   &  1.2943   &  0.7708   &  1.2943     \\
    \hline
    \hline
    $v_{1}$&$v_{2}$&$v_{3}$&$v_{4}$&$v_{5}$&$v_{6}$&$v_{10}$&$v_{11}$&$v_{12}$ \\
    \hline
    3.5308   &  2.1900  &  2.1900   & 3.5308 &   2.0210  &  2.0210  &  1.2943 &   0.7708 &   1.2943     \\
    \hline
    \end{tabular}
    \label{tab:my_label}
\end{table}

For the system (\ref{system}), we take its parameter $\lambda$ to increase from $1$ to $10^7$ and get the corresponding numerical solutions $(u_{\lambda}, v_\lambda)$. We find that with the increase of this parameter, $u_{\lambda}$ tends to $0$ at vertices in $\Omega_{a}^{c}$ and $v_{\lambda}$ tends to $0$ at vertices in $\Omega_{b}^{c}$, which are shown in Figure 2(a) and 2(b). In order to show the broken lines of the function values more clearly, we only select the values of several representative points to draw in the figure. Furthermore, the values of $u_\lambda$ in $\Omega_{a}$ and $v_\lambda$ in $\Omega_{b}$ just converge to those corresponding values in Table 1, which are shown in  Figure 2(c) and 2(d).

These above numerical results are completely consistent with our theorems proved in the previous sections.

\begin{figure}[htbp]
\begin{minipage}{0.5\textwidth}
  \centerline{\includegraphics[width=0.8\textwidth]{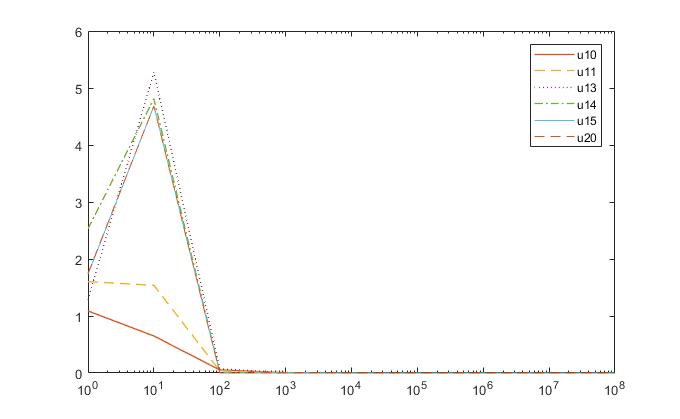}}
  \centerline{(a) Trend of $u_{\lambda}$ on $\Omega_{a}^{c}$}
\end{minipage}
\hfill
\begin{minipage}{0.5\textwidth}
  \centerline{\includegraphics[width=0.8\textwidth]{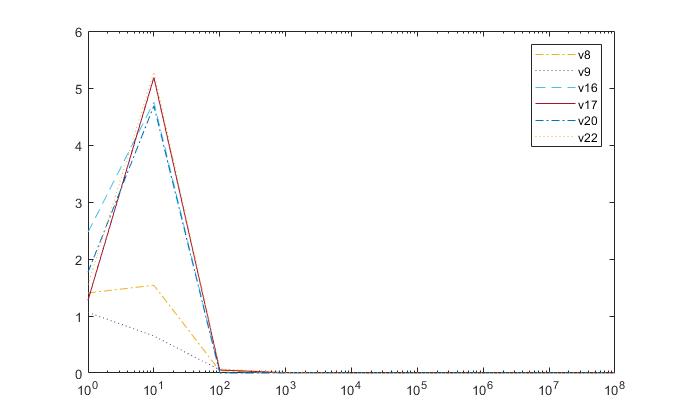}}
  \centerline{(b) Trend of $v_{\lambda}$ on $\Omega_{b}^{c}$}
\end{minipage}
\vfill
\begin{minipage}{0.5\textwidth}  
  \centerline{\includegraphics[width=0.8\textwidth]{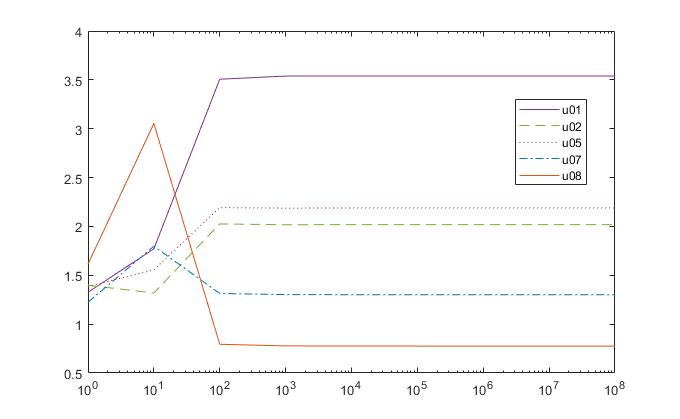}}
  \centerline{(c) Trend of $u_{\lambda}$ on $\Omega_{a}$}
\end{minipage}
\hfill
\begin{minipage}{0.5\textwidth}
  \centerline{\includegraphics[width=0.8\textwidth]{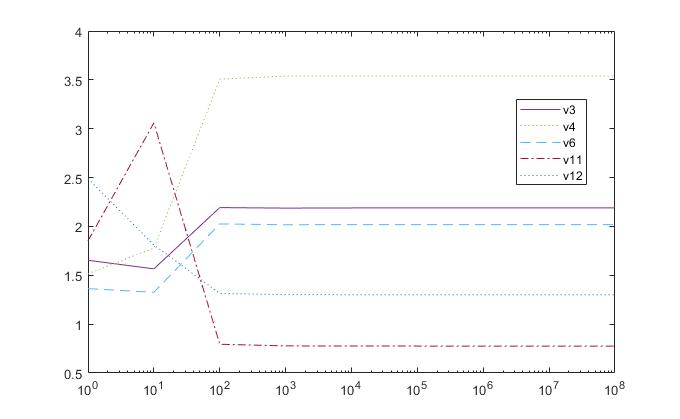}}
  \centerline{(d) Trend of $v_{\lambda}$ on $\Omega_{b}$}
\end{minipage}
\caption{Trends of the numerical solutions}
\label{fig:ob}
\end{figure}

\end{document}